\documentclass[12pt,leqno,twoside]{amsart}
\usepackage{amssymb}
\usepackage{amsmath}
\usepackage{latexsym}
\usepackage{verbatim}
\usepackage{epsfig}
\usepackage{xypic}

\newtheorem{theorem}{Theorem}[section]

\newtheorem{proposition}[theorem]{Proposition}
\theoremstyle{definition}
\newtheorem{definition}[theorem]{Definition}
\theoremstyle{remark}
\newtheorem{remark}[theorem]{\sc Remark}
\newtheorem{example}[theorem]{\sc Example}
\newcommand{\cO}{{\mathcal O}}

\newcommand{\ve}{{\varepsilon}}

\newcommand{\bR}{{\mathbb R}}
\newcommand{\bC}{{\mathbb C}}
\newcommand{\bP}{{\mathbb P}}

\newcommand{\bZ}{{\mathbb Z}}

\begin{document}

\title[Homological index for 1-forms and a Milnor number]
{Homological index for 1-forms and a Milnor number for
isolated singularities}

\author{\sc W. Ebeling}

\address{Wolfgang Ebeling: Universit\"at Hannover, Institut f\"ur Mathematik,
 Postfach 6009, D-30060 Hannover, Germany.}

\email{ebeling@math.uni-hannover.de}

\author{\sc S. M. Gusein-Zade}

\address{Sabir M. Gusein-Zade: Moscow State University, Faculty of Mechanics
 and Mathematics, Moscow, GSP-2, 119992, Russia.}

\email{sabir@mccme.ru}

\author{\sc J. Seade}

\address{Jos\'e Seade: Instituto de Matem\'aticas, Unidad Cuernavaca,
Universidad Nacional Aut\'onoma de M\'exico, Apartado postal 273-3,
C.P. 62210, Cuernavaca, Morelos, M\'exico.}

\email{jseade@matem.unam.mx}

\thanks{The research was partially supported by the DFG programme
"Global methods in complex geometry". The second author also was
partially supported by the grants RFBR--01--01--00739, INTAS--00--259;
the third author had partial support by CONACYT grant G36357/E.
}


%\keywords{}

%\date{April 10, 2003}

%\dedicatory{}

%\commby{}

%%% ----------------------------------------------------------------------
\begin{abstract} 
We introduce a notion of a homological index of a holomorphic 1-form
on a germ of a complex analytic variety with an isolated singularity,
inspired by X. G\'omez-Mont and G.-M. Greuel. For isolated complete
intersection singularities it coincides with the index defined earlier
by two of the authors. Subtracting from this index another one,
called radial, we get an invariant of the singularity which does not
depend on the 1-form. For isolated complete intersection singularities
this invariant coincides with the Milnor number. We compute this
invariant for arbitrary curve singularities and compare it with the
Milnor number introduced by R.-O. Buchweitz and G.-M. Greuel for
such singularities.
\end{abstract}
%
%%% ---------------------------------------------------------------------

\maketitle
%%% ----------------------------------------------------------------------

\setcounter{section}{0}
\section{Introduction}\label{intro}
%%%%%%%%%%%%%%%%%%%%%%%%%%%%%%

Indices of vector fields on singular varieties have been studied in
a number of papers, e.g., \cite {ASV, Go, KT, Sch, SS1, SS2}. One problem
with this approach is that the conditions for a vector field on the
ambient space to be tangent to a singular space are rather strong and
this makes their study difficult. On the contrary a $1$-form on the
ambient space always defines a 1-form on a subvariety. This and other
facts motivated the study of indices of $1$-forms on singular varieties
in \cite {EG2, EG3}.
 
Here we introduce a notion of a homological index of a 1-form on a
complex analytic variety with an isolated singular point analogous to
the homological index for vector fields defined in \cite {Go}. Due to
results of \cite {Gr}, the homological index of a 1-form on an isolated
complete intersection singularity is identified with the index defined
in \cite {EG2, EG3}. For a 1-form on an isolated complex analytic
singularity there is defined another index, which we call
radial. This index, as well as the homological one, satisfy the law of
conservation of number. For the homological index this follows from
\cite{GG}. For 1-forms on smooth varieties both indices are equal to the
usual Poincar\'e--Hopf index. This implies that the difference of these
two indices does not depend on the 1-form and is an invariant of the
singularity. For isolated complete intersection singularities we identify
it with the Milnor number. It is natural to try to compare the invariant
of an isolated singular point of a variety introduced in this paper with
other invariants coinciding with the Milnor number for isolated complete
intersection singularities. 

For curve singularities a notion of the Milnor number was introduced by
R.-O. Buchweitz and G.-M. Greuel in \cite{BG}. For a smoothable curve
singularity $(C,0)$ it is equal to $1-\chi(\widetilde C)$ for a smoothing
$\widetilde C$ of the singularity, $\chi(\cdot)$ is the Euler
characteristic. However, there exist surface singularities which have
smoothings with different Euler characteristics (\cite{P}). This does not
permit to generalize the notion of the Milnor number to higher dimensions
so that for smoothable singularities it has the usual expression in terms
of the Euler characteristic. We compute our invariant for curve
singularities and compare it with the Milnor number of R.-O. Buchweitz and
G.-M. Greuel.

The authors are thankful to G.-M. Greuel for useful discussions
and to H.-C. von Bothmer for computations using {\sc Macaulay}. The last
two authors are grateful for the hospitality of the University of
Hannover where this work was carried out.

\bigskip
%%%%%%%%%%%%%%%%%%%%%%%%%%%%%% %%%%%%%%%%%%%%%%%%%%%%%%%%%%%%
\section{Indices of $1$-forms}
%%%%%%%%%%%%%%%%%%%%%%%%%%%%%% %%%%%%%%%%%%%%%%%%%%%%%%%%%%%%

In \cite {EG3} there was defined a notion of an index of a 1-form on a (real) manifold
with isolated singularities. Here we shall consider complex analytic varieties with
isolated singularities and (complex) 1-forms on them. Therefore we reformulate the definition
from \cite {EG3} in this setting.

Let $(V,0)\subset (\bC^N, 0)$ be a germ of a purely $n$-dimensional complex
analytic variety with an isolated singularity at the origin $0$.
Let $\omega$ be a 1-form on $V$ with an isolated singularity at the origin $0$.
Here this means that $\omega$ is a continuous, nowhere-vanishing section of
the complex cotangent bundle of $V-\{0\}$.
We consider an index that measures the "lack of radiality" of such a $1$-form.
This is similar to the index for vector fields considered in \cite{ASV, EG1, KT, SS2}.

Let us fix a radial vector field $v_{\rm rad}$ on $(V,0)$, e.g., the gradient on the
smooth part of $V$ of the real valued function $\Vert z\Vert$ with respect to a
Riemannian metric. A radial 1-form $\omega_{\rm rad}$ is a (complex) 1-form on $V$ whose
value on the radial vector field $v_{\rm rad}$ has positive real part at each point
in a punctured neighbourhood of the origin $0$ in $V$. The space of such 1-forms
is connected.

Let $\omega_1$ and $\omega_2$ be 1-forms on $(V,0)$ with isolated
singularities at the origin. Choose $\ve > \ve' > 0$ sufficiently small,
let $K_{\ve} = V \cap S_{\ve}$ and $K_{\ve'} = V \cap S_{\ve'}$ be
the corresponding links, and let $Z$ be the cylinder $V \cap [B_{\ve}
\setminus {\rm Int}(B_{\ve'})]$, where $B_\rho$ is the ball of radius
$\rho$ around the origin $0$ in $\bC^N$, $S_\rho$ is its boundary.
Let $\widetilde\omega$ be a 1-form on the cylinder $Z$ which coincides
with $\omega_1$ in a neighbourhood of $K_{\ve}$ and with $\omega_2$ in
a neighbourhood of $K_{\ve'}$ and which
has isolated singular points $Q_1$, ..., $Q_s$ inside $C$. The sum
$d(\omega_1,\omega_2)$ of the (usual) local indices
${\rm ind}_{Q_i}\, \widetilde\omega$ of the form $\widetilde\omega$
at these points depends only on the forms $\omega_1$ and $\omega_2$
and will be called {\em the difference} of these forms. One has
$d(\omega_1,\omega_2) = -d(\omega_2,\omega_1)$.

\begin{definition}\label{radial index}
The {\em radial index} ${\rm Ind}_{\rm rad}(\omega;V,0)$ (or simply ${\rm
Ind}_{\rm rad}\, \omega$) of the 1-form
$\omega$ on $V$ is defined by the equation
$${\rm Ind}_{\rm rad}\, \omega \,=\, (-1)^n + d(\omega,\omega_{\rm rad})\,.$$ 
\end{definition}

\begin{remark}
%% {\bf 1.}
One can see that the index of the radial 1-form
$\omega_{\rm rad}$ is equal to $(-1)^n$. The sign is chosen so that this index
coincides with the usual one when $V$ is smooth at $0$.
\end{remark}

\begin{remark}
%% \noindent{\bf 2.}
There is a one-to-one correspondence between complex 1-forms
on a complex analytic manifold $V-\{0\}$ and real 1-forms on it. Namely,
to a complex 1-form $\omega$ one associates the real 1-form
$\eta={\rm Re}\,\omega$; the 1-form $\omega$ can be restored from $\eta$
by the formula $\omega(v)=\eta(v)-i\eta(iv)$ for $v\in T_x(V-\{0\})$.
This explains why the radial index of a complex 1-form can be expressed
through the corresponding index of its real part defined in \cite {EG3}.
In other words the radial index ${\rm Ind}_{\rm rad}\, \omega$ of a complex
1-form $\omega$ can be defined as $(-1)^n$-times the radial index of
the real part of $\omega$.
\end{remark}

\begin{remark}
The radial index obviously satisfies the law of conservation of number,
i.e., if $\omega'$ is a 1-form on $V$ close to $\omega$, then: 
$$ {\rm Ind}_{\rm rad}(\omega; V,0) = {\rm Ind}_{\rm rad}(\omega'; V,0) + \sum
{\rm Ind}_{\rm rad}(\omega';V,x) \,,$$ where the sum on the right hand
side is over all those points $x$ in a small punctured neighbourhood
of the origin $0$ in $V$ where the form $\omega'$ vanishes
(this follows from the fact that $d(\omega_1, \omega_3)=
d(\omega_1, \omega_2)+d(\omega_2, \omega_3)$).
\end{remark}

\begin{example}\label{ex}
Let $\omega$ be a holomorphic 1-form on a curve singularity $(C,0)$ with
$C=\cup_{i=1}^r C_i$, where $C_i$ are the irreducible components of $C$.
Let $t_i$ be a uniformization parameter on the component $C_i$ and
let the restriction $\omega_{\vert C_i}$ be of the form $(a_it_i^{m_i} +\
terms\ of\ higher\ degree)dt_i$, $a_i\ne 0$. Then ${\rm Ind}_{\rm rad}\,
\omega_{\vert C_i}=m_i$. Therefore $d(\omega_{\vert C_i}, \omega_{\rm
rad}{}_{\vert C_i})=m_i+1$,
$d(\omega, \omega_{\rm rad})=\sum_{i=1}^r (m_i+1)$, ${\rm Ind}_{\rm rad}\,
\omega=
\sum_{i=1}^r m_i + (r-1)$.
\end{example}

\begin{example}\label{ex2}
Let $(W,0)\subset(\bR^M,0)$ be a germ of a (real) analytic
variety with an isolated singular point at the origin and
let $f:(W,0)\to(R,0)$ be a germ of a real analytic function
on $(W,0)$ with an isolated critical point at the origin
(the last means that the differential $df$ has an isolated
singular point at the origin on $V$). One can show that
the (radial) index of the differential $df$ (as defined in
\cite{EG3} for real 1-forms) is equal to $1-\chi(F_-)$, where
$F_-$ is the "negative" Milnor fibre
$\{f=-\delta\}\cap B_\varepsilon$ for
$0<\delta\ll\varepsilon$ small enough. This is an analogue for
singular varieties of a formula of \cite{Ar}. 
Let $f$ be a germ
of a holomorphic function on a complex analytic isolated
singularity $(V,0)$ with an isolated critical point at the
origin. The map $f$ defines a Milnor fibration (as noticed, e.g., by
H.~Hamm \cite{Ha}). One can show that the Milnor fibre of $f$ is homotopy
equivalent to the Milnor fibre of its real part ${\rm Re}\, f$. Therefore
the radial index
${\rm Ind}_{\rm rad}\, df$ of its differential $df$ is equal to
$(-1)^n(1-\chi(F))$, where
$F$ is the Milnor fibre $\{f=\delta\}\cap B_\varepsilon$
of the function $f$ on $V$
($0<\vert\delta\vert\ll\varepsilon$ are small enough).
\end{example}

Now let $(V,0)\subset (\bC^{n+k},0)$ be an isolated complete
intersection singularity, defined by an analytic map
$$f = (f_1,...,f_k) : (\bC^{n+k},0) \to (\bC^k,0) \,.$$ 
For a 1-form $\omega$ on $V$ with an isolated singularity at the origin $0$, the
covectors $\omega$, $df_1$, ..., $df_k$ are linearly independent at each point of $V$
away from $0$. This defines a map $(\omega, df_1,...,df_k)$ from 
$V \setminus \{0\}$ into the Stiefel manifold $W^*_{k+1}(\bC^{n+k})$ of $(k+1)$-frames
in the dual $\bC^{n+k}$. The manifold $W^*_{k+1}(\bC^{n+k})$ is $(2n-2)$-connected
and its first non-zero homotopy group $\pi_{2n-1}(W^*_{k+1}(\bC^{n+k})$ is
isomorphic to $\bZ$ (see, e.g., \cite {Hu}). Therefore the restriction of this map to
the link $K_\varepsilon$ of $(V,0)$ has a degree: that of the map induced in the homology of dimension $(2n-1)= \dim K_\varepsilon$.

\begin{definition} The {\em index} ${\rm Ind}(\omega;V,0)=
{\rm Ind}\,\omega$ of the 1-form $\omega$ is the degree of the map 
$$(\omega, df_1,...,df_k) : K_\varepsilon \to W^*_{k+1}(\bC^{n+k}) \,.$$
\end{definition}

This index, defined in \cite{EG2, EG3}, is an analogue of the
GSV-index for vector fields, introduced in \cite {GSV,SS1}.
It is proved in \cite{EG3}, Proposition 3, that this index
equals the number of zeroes, counted with multiplicities,
of any extension of $\omega$ to a Milnor fibre
$V_t = f^{-1}(t) \cap B_{\ve}$ of $(V,0)$. Let $\mu (V,0) $
be the Milnor number of the isolated complete intersection
singularity $(V,0)$.

\begin{proposition} \label{p: Milnor Number}
For a $1$-form $\omega$ on $(V, 0)$ with an isolated 
singularity at the origin $0$ one has
$${\rm Ind}\,\omega = \mu (V,0) + d(\omega,\omega_{\rm rad}) + (-1)^n$$
and therefore the difference ${\rm Ind}\,\omega - {\rm Ind}_{\rm rad}\, \omega$
between its index and its radial index is independent of $\omega$ and is equal
to the Milnor number $\mu (V,0)$.
\end{proposition}

An analogue of this proposition for vector fields also holds:
\cite{SS2}, Proposition 1.4. If the 1-form $\omega$ is holomorphic, i.e.,
if it is the restriction to $(V, 0)$ of a holomorphic 1-form
$\omega = \sum_{i=1}^{n+k}A_i(x)dx_i$ on $(\bC^{n+k}, 0)$, this index
has an algebraic expression as the dimension of a certain algebra
\cite{EG2, EG3}. Namely, 
$${\rm Ind}\,\omega=\dim_\bC \cO_{\bC^N,0} / I\,,$$
where $I$ is the ideal generated by $f_1$, ..., $f_k$
and the $(k+1) \times (k+1)$-minors of the matrix:
$$\begin{pmatrix} \frac{\partial f_1}{\partial x_1} & {\cdots } &
 \frac{\partial f_1}{\partial x_{n+k}} \\ {\vdots} & {\cdots} & 
{\vdots} \\ \frac{\partial 
f_k}{\partial x_1} & \cdots &
 \frac{\partial f_k}{\partial x_{n+k}} \\
A_1 & \cdots & A_{n+k}
\end{pmatrix}\,.
$$
This formula was obtained by L\^e D.T. and G.-M. Greuel for the case when $\omega$ is the
differential of a function (\cite {Gr, Le}).
\bigskip

%%%%%%%%%%%%%%%%%%%%%%%%%%%%%% %%%%%%%%%%%%%%%%%%%%%%%%%%%%%%
\section{The homological index}
%%%%%%%%%%%%%%%%%%%%%%%%%%%%%% %%%%%%%%%%%%%%%%%%%%%%%%%%%%%%

Let $(V,0)\subset (\bC^N, 0)$ be any germ of an analytic variety of pure
dimension $n$ with an isolated singular point at the origin (not
necessarily a complete intersection). Given a holomorphic form
$\omega$ on $(V,0)$ with an isolated singularity, we consider
the complex $(\Omega^\bullet_{V,0}, \wedge\omega)$: 
$$
0 \to \cO_{V,0}
\to \Omega^1_{V,0} \to ... \to \Omega^n_{V,0} \to 0\,, 
$$
where $\Omega^i_{V,0}$ are the modules of sheaves of differential forms on $(V,0)$
and the arrows are given by the exterior product by the form $\omega$.

This complex is the dual of the Koszul complex considered in \cite{Go}, (1.4).
It was used by G.M.~Greuel in \cite{Gr} for complete intersections.
The sheaves $\Omega^i_{V,0}$ are coherent sheaves and the homology groups
of the complex $(\Omega^\bullet_{V,0}, \wedge\omega)$ are concentrated at
the origin and therefore are finite dimensional. The definition below 
was inspired by that of \cite{Go} for vector fields.

\begin{definition} The {\em homological index}
$\,{\rm Ind}_{\rm hom}(\omega; V,0)={\rm Ind}_{\rm hom}\, \omega$ of the 1-form
$\omega$ on $(V, 0)$ is $(-1)^n$ times the Euler characteristic of the above
complex:
\begin{equation}\label{eq1}
{\rm Ind}_{\rm hom}(\omega; V,0) ={\rm Ind}_{\rm hom}\, \omega= \sum_{i=0}^n
(-1)^{n-i} h_i(\Omega^\bullet_{V,0},\wedge\omega)\,,
\end{equation}
where $h_i(\Omega^\bullet_{V,0},\wedge\omega)$ is the dimension of the
corresponding homology group as a vector space over $\bC$.
\end{definition}

\begin{theorem} \label{t:homological} Let $\omega$ be a holomorphic 
1-form on $V$ with an isolated singularity at the origin $0$. 

\medskip
\item {\rm(i)} If $V$ is smooth, then ${\rm Ind}_{\rm hom}\, \omega$
equals the usual local index of the 1-form $\omega$.

\medskip
\item {\rm(ii)} The homological index satisfies the law of conservation of number: if
$\omega'$ is a holomorphic 1-form on $V$ close to $\omega$, then: 
$${\rm Ind}_{\rm hom}(\omega; V,0) = {\rm Ind}_{\rm hom}(\omega'; V,0)
 + \sum {\rm Ind}_{\rm hom}(\omega';V,x) \,,$$ 
where the sum on the right hand side is over all those points $x$ 
in a small punctured neighbourhood of the origin $0$ in $V$ where the form $\omega'$ vanishes.

\medskip
\item {\rm(iii)} If $(V,0)$ is an isolated complete intersection
singularity, then the homological index ${\rm Ind}_{\rm hom}\, \omega$ coincides
with the index ${\rm Ind}\,\omega$.
\end{theorem}

\begin{proof}
Statement (i) is straightforward and it is a special case of
statement (iii). Statement (ii) is a particular case of the
main theorem in \cite{GG}. For statement (iii), we notice
that (on an isolated complete intersection singularity $(V,0)$)
the index ${\rm Ind}\,\omega$ also satisfies the law of
conservation of number and coincides with the homological
index ${\rm Ind}_{\rm hom}\,\omega$ on smooth varieties. This
implies that the difference between these two indices is a
locally constant, and therefore constant, function on the
space of 1-forms on $(V, 0)$ with an isolated singular point
at the origin. 
Therefore it suffices to prove (iii) for $\omega = df$ where $f$ is a
holomorphic function on $(V,0)$ with an isolated critical point at the
origin. The de Rham lemma in \cite{Gr}, Lemma 1.6, implies that the
homology groups of the complex
$(\Omega^\bullet_{V,0}, \wedge df)$ vanish in dimensions
$i = 0, 1, ..., n-1$. The statement then follows from the
Remark after Lemma 5.3 of \cite{Gr} (see also \cite{EG3}). 
\end{proof}

\begin{remark}
The minimal value of the homological index ${\rm Ind}_{\rm hom}\,\omega$ is
attained by restrictions to $V$ of generic 1-forms on $\bC^N$
which do not vanish at the origin. The subset of forms with
this index in $\Omega^1_{V,0}$ is open, dense and connected.
Moreover, each 1-form $\omega$ can be approximated by a 1-form,
the index of which at the origin coincides with the minimal one
and all its zeroes on $V-\{0\}$ are non-degenerate. This approximation
can be chosen of the form $\omega + \ve d\ell$ for a linear function
$\ell$.
\end{remark}

\begin{remark} 
We notice that one has 
an invariant for functions on $(V, 0)$ with an isolated singularity at the
origin defined by $f \mapsto {\rm Ind}_{\rm hom}\,df$. By the theorem above, if 
$(V,0)$ is an isolated complete intersection singularity, this invariant
%% measures
counts the number of critical points of the function $f$ on a Milnor fibre.
\end{remark}

\begin{remark}\label{remMvS}
Let $(C,0)$ be a curve singularity and let $(\bar{C},\bar{0})$
be its normalization. Let $\tau=\dim \rm {Ker}(\Omega^1_{C,0}\to
\Omega^1_{\bar{C},\bar{0}})$, $\lambda=\dim \omega_{C,0}/
c(\Omega^1_{C,0})$, where $\omega_{C,0}$ is the dualizing
module of Grothendieck, $c:\Omega^1_{C,0}\to\omega_{C,0}$
is the class map (see \cite{BG}). In \cite{MS} there is considered
a Milnor number of a function $f$ on a curve singularity
introduced by V.Goryunov. One
can see that this Milnor number can be defined for a 1-form
$\omega$ with an isolated singularity on $(C,0)$ as
well (as $\dim \omega_{C,0}/\omega\wedge \cO_{C,0}$) and is
equal to ${\rm Ind}_{\rm hom}\,\omega + \lambda - \tau$.
\end{remark}
\bigskip

%%%%%%%%%%%%%%%%%%%%%%%%%%%%%% %%%%%%%%%%%%%%%%%%%%%%%%%%%%%%
\section{A generalized Milnor number}
%%%%%%%%%%%%%%%%%%%%%%%%%%%%%% %%%%%%%%%%%%%%%%%%%%%%%%%%%%%%

We recall that, by Proposition \ref{p: Milnor Number}, for $1$-forms on 
an isolated complete intersection singularity the difference between its
index and its radial index is the Milnor number of the singularity,
independently of the choice of the $1$-form. On a germ $(V,0)$ in general,
the index ${\rm Ind}\,\omega$ of a 1-form is not defined, but the 
homological index ${\rm Ind}_{\rm hom}\,\omega$ is and it coincides with the
index ${\rm Ind}\,\omega$ for complete intersection germs. The radial
index is always defined.

The laws of conservation of numbers for the homological and the radial
indices of 1-forms together with the fact that these two indices
coincide on smooth varieties imply that their difference is a
locally constant, and therefore constant, function on the space
of 1-forms on $V$ with isolated singularities at the origin.
Therefore one has the following statement.

\begin{proposition} Let $(V,0)$ be a germ of a complex analytic
space of pure dimension $n$ with an isolated singular point at
the origin. Then the difference
$${\rm Ind}_{\rm hom}\,\omega - {\rm Ind}_{\rm rad}\,\omega$$
between the homological and the radial indices
does not depend on the 1-form $\omega$.
\end{proposition}

This proposition, together with \ref{p: Milnor Number},
permits to consider the difference 
$$\nu(V,0) = {\rm Ind}_{\rm hom}(\omega;V,0) - {\rm Ind}_{\rm rad}(\omega;V,0)$$
as a generalized Milnor number of the singularity $(V,0)$.

\begin{remark}\label{diff}
One can define the invariant $\nu(V,0)$ using any 1-form
$\omega$ on $V$ with an isolated singular point at the
origin, say, the differential $df$ of a holomorphic function
$f$ on $V$. In this case one can describe $\nu(V,0)$ in
somewhat different terms. When $(V,0)$ is an isolated complete
intersection singularity, the homological index ${\rm
Ind}_{\rm hom}\,df$ is equal to the sum $\mu(V,0) + \mu(f)$ where
$\mu(f)$ is the Milnor number of the germ $f$ on $(V,0)$. Thus
for an arbitrary isolated singularity $(V, 0)$ and for a
function germ $f$ on it, the homological index
${\rm Ind}_{\rm hom}\,df$ can be considered as a generalization of
this sum. 
The Milnor number $\mu(f)$ of a germ $f$ of a holomorphic function 
with an isolated critical point
can be defined for any isolated singularity
$(V, 0)$ of pure dimension $n$ as $(-1)^{n}(1-\chi(F))$,
where $F$ is the Milnor fibre
of $f$, and is equal to the radial index ${\rm Ind}_{\rm rad}\,df$ 
(see Example~\ref{ex2}). Therefore the invariant $\nu(V,0)$ is a
generalization of the number $\mu(V,0) + \mu(f)$ minus the Milnor number
$\mu(f)$.
\end{remark}

There are other invariants of isolated singularities of complex analytic
varieties which coincide with the Milnor number for isolated complete
intersection singularities. One of them is $(-1)^n$ times the reduced
Euler characteristic (i.e., the Euler characteristic minus 1) of the
absolute de Rham complex of $(V, 0)$. It is natural to try to compare
$\nu(V,0)$ with them.

\begin{theorem}
For a curve singularity $(C,0)$,
$$
\nu(C, 0)=\dim_\bC \Omega^1_{C,0}/d\cO_{C,0}\,.
$$
\end{theorem}

\begin{proof}
Let $\pi:(\bar{C}, \bar{0}) \to (C,0)$ be the normalization of
the curve $(C, 0)$, $\bar{0}=\pi^{-1}(0)$. Let $C=\cup_{i=1}^r C_i$
be the decomposition of the curve $(C,0)$ into the union of
irreducible components, let $t_i$ be the uniformization
parameter on $C_i$, and let $\omega_{\vert C_i} =
(a_it_i^{m_i}+\ terms\ of\ higher\ degree)dt_i$, $a_i\ne 0$.
Consider two commutative diagrams with exact rows
$$\diagram
0 \rto & \cO_{C,0} \rto^{ \wedge \omega}
\dto_{\pi_0^\ast} & \Omega^1_{C,0} \rto \dto_{\pi_1^\ast}
& \Omega^1_{C,0}/\omega \wedge \cO_{C,0} \rto \dto & 0 \\
0 \rto & \cO_{\bar{C}, \bar{0}} \rto^{\wedge \bar{\omega}} &
\Omega^1_{\bar{C}, \bar{0}} \rto &
\Omega^1_{\bar{C},\bar{0}}/\bar{\omega} \wedge
\cO_{\bar{C},\bar{0}} \rto & 0
\enddiagram
$$
$$\diagram
0 \rto & {\mathfrak m}_{C,0} \rto^d
\dto_{\widetilde{\pi}_0^\ast} & \Omega^1_{C,0} \rto \dto_{\pi_1^\ast}
& \Omega^1_{C,0}/d \cO_{C,0} \rto \dto & 0 \\
0 \rto & {\mathfrak m}_{\bar{C}, \bar{0}} \rto^d &
\Omega^1_{\bar{C}, \bar{0}} \rto &
0
%% \Omega^1_{\bar{C},\bar{0}}/d \cO_{\bar{C},\bar{0}}
 \rto & 0
\enddiagram
$$
where ${\mathfrak m}_{C,0}$ is the maximal ideal in the ring
$\cO_{C,0}$, ${\mathfrak m}_{\bar C,\bar 0}$ is the ideal of
germs of functions on the normalization $(\bar{C}, \bar{0})$
equal to zero at all the point from $\bar 0$. Note that
the homomorphisms $\pi_0^*$ and $\widetilde\pi_0^*$ are
injective. Here $\dim \Omega_{C,0}^1/\omega\wedge\cO_{C,0}$
is equal to the homological index ${\rm Ind}_{\rm hom}\,\omega$,
$\dim \Omega_{\bar C,\bar 0}^1/\bar\omega\wedge\cO_{\bar C,
\bar 0}=\sum m_i={\rm Ind}_{\rm rad}\,\omega - (r-1)$ (see Example
\ref{ex}), $\dim {\rm coker}\,\pi^*_0 = \dim
{\rm coker}\,\widetilde{\pi}^*_0 + (r-1)$. Applying the Snake Lemma (see,
e.g., \cite{Ei}) to these diagrams (all kernels and cokernels of the
vertical homomorphisms are finite dimensional) we get the statement.
\end{proof}

\begin{remark} 
A notion of a generalized Milnor number of a curve singularity $(C,0)$
was introduced in \cite{BG} as $\dim_\bC \omega_{C,0}/d\cO_{C,0}$,
where $\omega_{C,0}$ is the dualizing module of Grothendieck. For
smoothable curve singularities, it is equal to $1-\chi(\widetilde C)$,
where $\widetilde C$ is a smoothing of $(C,0)$. (All smoothings
of a curve singularity have the same Euler characteristic.) 
From the proof of \cite{BG}, Theorem 6.1.3, it follows that the
Milnor number defined by R.-O. Buchweitz and G.-M. Greuel is
equal to $\nu(C,0) + \lambda - \tau$, where $\tau$ and
$\lambda$ are defined in Remark \ref{remMvS}. For
complete intersection curve singularities $ \lambda = \tau$.
\end{remark}

\begin{example}
Let $(S,0)\subset(\bC^5,0)$ be the cone over the rational normal curve
in $\bC\bP^4$ (this is Pinkham's example \cite{P} mentioned in the
Introduction). According to Example \ref{ex2}, for a generic linear
function $\ell$ on $\bC^5$, the radial index ${\rm Ind}_{\rm rad}\,d\ell$
of its differential on the surface $S$ is equal to the Milnor number
$\mu(\ell_{\vert S})$ (see Remark \ref{diff} for the definition of the
latter one). This Milnor number can be easily computed (say, using the
formula of N. A'Campo \cite{AC}) and is equal to $3$. All the equations of
the surface $S$ are homogeneous of degree $2$ and therefore the modules
$\Omega_{S,0}^i$,
$i=0,1,2$, have natural gradings (where we consider the
differentials $dx_j$ of the variables having degree $0$,
$j=1, \ldots, 5$). The differentials in the complex
$(\Omega^\bullet_{S,0}, \wedge d\ell)$ respect the grading. The
differentials in the absolute de Rham complex $(\Omega^\bullet_{S,0}, d)$
have degree $(-1)$. Therefore the Euler characteristics of these
complexes can be computed from the Poincar\'e series $P_i(t)$ of the
gradings on the modules $\Omega_{S,0}^i$. These Poincar\'e series
(computed using {\sc Macaulay}) are equal to:
\begin{eqnarray*}
P_0(t) & = &
1+5t+9t^2+13t^3+17t^4+21t^5+25t^6+29t^7+33t^8+37t^9+ \ldots \, ,
\\
P_1(t) & = &
5+19t+24t^2+32t^3+40t^4+48t^5+56t^6+64t^7+72t^8+80t^9+ \ldots \, ,
\\
P_2(t) & = &
10+21t+15t^2+19t^3+23t^4+27t^5+31t^6+35t^7+39t^8+43t^9+ \ldots \, .
\end{eqnarray*}
Therefore the Euler characteristic of the complex $(\Omega^\bullet_{S,0},
\wedge d\ell)$ (equal to the homological index ${\rm Ind}_{\rm hom}\,d\ell$)
is equal to $\left(P_0(t)-P_1(t)+P_2(t)\right)_{\vert t\mapsto 1}=13$ and
the invariant
$\nu(S,0)$ introduced above is equal to $10$. One can see that the reduced
Euler characteristic of the absolute de Rham complex of $(S,0)$ is
equal to $\left(P_0(t)-tP_1(t)+t^2P_2(t)\right)_{\vert t\mapsto 1}$
and hence also equal to $10$. 
\end{example}

\begin{remark}
Note that the differentials of the complex
$(\Omega^\bullet_{V,0}, \wedge\omega)$ are homomorphisms of
$\cO_{V,0}$--modules while those of the de Rham complex do not have this
property.
\end{remark}


\bigskip
\begin{thebibliography}{MM}

\bibitem[1]{AC} {\em N. A'Campo,} La fonction z\^eta d'une monodromie.
Comment. Math. Helv. {\bf 50} (1975), 233--248.

\bibitem[2]{ASV} {\em M. Aguilar, J. Seade, A. Verjovsky,} Indices of
vector fields and topological invariants or real analytic singularities.
J. Reine Angew. Math. {bf 504} (1998), 159--176.

\bibitem[3]{Ar} {\em V.I. Arnold,}
Index of a singular point of a vector field, the Petrovskii--Oleinik inequality,
and mixed Hodge structures. Funct. Anal. Appl. {\bf 12} (1978), no.1, 1--12. 

\bibitem[4]{BG} {\em R.-O. Buchweitz, G.-M. Greuel,} The Milnor number
and deformations of complex curve singularities. Invent. Math. {\bf 58}
(1980), 241--281. 

\bibitem[5]{EG1} 
{\em W. Ebeling, S. M. Gusein-Zade,} On the index of 
a vector field at an isolated singularity. In: The Arnoldfest, edited by
E. Bierstone et al., Fields Inst. Commun. {\bf 24} (1999), 141--152, AMS.

\bibitem[6]{EG2} {\em W. Ebeling, S. M. Gusein-Zade,} On the index of a 
holomorphic $1$-form on an isolated complete intersection singularity.
Doklady Math. {\bf 64} (2001), 221--224.

\bibitem[7]{EG3} {\em W. Ebeling, S. M. Gusein-Zade,} Indices of 1-forms 
on an isolated complete intersection singularity. To appear in Moscow
Math. Journal.

\bibitem[8]{Ei}
{\em D. Eisenbud,} Commutative Algebra with a view toward
Algebraic Geometry. Springer Verlag, Graduate texts in Math. {\bf 150} 
(1994).

\bibitem[9]{GG}
{\em L. Giraldo, X. G\'omez-Mont,} A law of conservation of number for
local Euler characteristics. Contemp. Math. {\bf 311} (2002), 251--259.

\bibitem[10]{Go}
{\em X. G\'omez-Mont,} An algebraic formula for the index 
of a vector field on a hypersurface with an isolated singularity.
J. Algebraic Geom. {\bf 7} (1998), 731--752. 

\bibitem[11]{GSV}
{\em X. G\'omez-Mont, J. Seade, A. Verjovsky,} The index of a
holomorphic flow with an isolated singularity. Math. Ann. {\bf 291}
(1991), 737--751.

\bibitem[12]{Gr}
{\em G.-M. Greuel,} Der Gau\ss-Manin-Zusammenhang 
isolierter Singularit\"aten von vollst\"andigen Durchschnitten.
Math. Ann. {\bf 214} (1975), 235--266.

\bibitem[13]{Ha}
{\em H. Hamm,} Lokale topologische Eigenschaften komplexer R\"aume. Math.
Ann. {\bf 191} (1971), 235--252. 

\bibitem[14]{Hu}
{\em D. Husemoller,} Fibre bundles. 2nd edition, Graduate Texts in 
Maths. {\bf 20}, Springer Verlag, 1975.

\bibitem[15]{KT} {\em H. King, D. Trotman,}
Poincar\'e-Hopf theorems on stratified sets. Preprint 1996.

\bibitem[16]{Le} {\em L\^ e D.T.,} Computation of the Milnor number of an
isolated singularity of a complete intersection. Funct. Anal. Appl. {\bf
8} (1974), 45--49.

\bibitem[17]{MS}
{\em D. Mond, D. Van Straten,} Milnor number equals Tjurina number
for functions on space curves. J. London Math. Soc. {\bf 63} (2001),
177--187.

\bibitem[18]{P}
{\em H.C. Pinkham,} Deformations of algebraic varieties with
$G_m$ action. Ast\'erisque {\bf 20} (1974).

\bibitem[19]{Sch}
{\em M.-H. Schwartz,} Classes caract\'eristiques d\'efinies par une
stratification d'une vari\'et\'e analytique complexe.
C.R. Acad. Sci. Paris {\bf 260} (1965), 3262--3264, 3535--3537.

\bibitem[20]{SS1} {\em J. Seade, T. Suwa,} A residue formula for the
index of a holomorphic flow. Math. Ann. {\bf 304} (1994), 345-360.

\bibitem[21]{SS2} {\em J. Seade, T. Suwa,} An adjunction formula for local
complete intersections. Internat. J. Math. {\bf 9} (1998), 759--768.

\end{thebibliography}
\end{document}